\documentclass[10pt]{amsart}
\usepackage{palatino,amssymb,mathrsfs}
\usepackage[hypertex]{hyperref}
\usepackage[author-year,sorted]{amsrefs}%
   \title[Capacities in Wiener Space]{Capacities
   in Wiener Space,\\
   Quasi-Sure Lower Functions,\\
   and Kolmogorov's {\boldmath$\e$}-Entropy}

   \author[D. Khoshnevisan]{Davar Khoshnevisan}
   \thanks{The research of D.~Kh.~ is partially supported
   by a grant from the NSF}
   \address{Department\@ of Mathematics\\The University\@ of Utah\\
      155 S.\@ 1400 E.\\Salt Lake City, UT 84112--0090}
   \email{davar@math.utah.edu}
   \urladdr{http://www.math.utah.edu/\~{}davar}
   \author[D.\@ A.\@ Levin]{David A.\@ Levin}
   \address{Department\@ of Mathematics\\The University\@ of Utah\\
      155 S.\@ 1400 E.\\Salt Lake City, UT 84112--0090}
   \email{levin@math.utah.edu}
   \urladdr{http://www.math.utah.edu/\~{}levin}
   \author[P.\@ M\'endez]{Pedro J. M\'endez-Hern\'andez} 
   \address{Department\@ of Mathematics\\The University\@ of Utah\\ 
      155 S.\@ 1400 E.\\Salt Lake City, UT 84112--0090}
   \curraddr{Escuela de Matem\'atica\\Universidad de Costa Rica\\
      San Pedro de Montes de Oca\\Costa Rica}
   \email{mendez@math.utah.edu,pmendez@emate.ucr.ac.cr}
   \urladdr{http://www.math.utah.edu/\~{}mendez}

\theoremstyle{plain}{
\newtheorem{theorem}{Theorem}[section]}
\theoremstyle{plain}{
\newtheorem{proposition}[theorem]{Proposition}}
\theoremstyle{plain}{
   \newtheorem{lemma}[theorem]{Lemma}}
\theoremstyle{plain}{
   }
\theoremstyle{plain}{
   \newtheorem{corollary}[theorem]{Corollary}}
\theoremstyle{definition}{
   \newtheorem{definition}[theorem]{Definition}}
\theoremstyle{definition}{
   }
\theoremstyle{remark}{
   \newtheorem{remark}[theorem]{Remark}}
\numberwithin{equation}{section}
\newcommand{\K}{\mathrm{K}}
\newcommand{\C}{\mathrm{cap}}
\newcommand{\M}{\mathrm{M}}

\newcommand{\Kdim}{\overline{\dim}_{_{\mathscr M}}}

\newcommand{\Pdim}{\dim_{_{\mathscr P}}}
\newcommand{\F}{\mathscr{F}}
\newcommand{\B}{\mathscr{B}}

\newcommand{\e}{\varepsilon}
\newcommand{\s}{\sigma}

\renewcommand{\P}{\mathrm{P}}
\newcommand{\ee}{\mathscr{e}}
\newcommand{\E}{\mathrm{E}}
\newcommand{\R}{\mathbf{R}}
\renewcommand{\ee}{\mathbf{e}}
\date{September 25, 2004}
\begin{document}
\begin{abstract}
   We propose a set-indexed family of capacities
   $\{\C_G\}_{G\subseteq\R_+}$
   on the classical Wiener space $C(\R_+)$.
   This family interpolates between the Wiener measure ($\C_{\{0\}}$)
   on $C(\R_+)$ and the standard capacity ($\C_{\R_+}$)
   on Wiener space.
   We then apply our capacities to characterize all
   quasi-sure lower functions in $C(\R_+)$.
   In order to do this we derive the following
   capacity estimate (Theorem~\ref{main-EST})
   which may be of independent interest:
   There exists a constant $a>1$
   such that for all $r>0$,
   \[
      \frac 1a \K_G(r^6)e^{-\pi^2/(8r^2)} \le
      \C_G \{f^*\le r\}
      \le a\K_G(r^6)e^{-\pi^2/(8r^2)}.
   \]
   Here, $\K_G$ denotes the Kolmogorov $\e$-entropy
   of $G$, and $f^*:=\sup_{[0,1]}|f|$.
\end{abstract}
\maketitle
%
\section{Introduction}

Let $C(\R_+)$ denote the collection of all
continuous  functions $f:\R_+\to\R$. We endow
$C(\R_+)$ with its usual topology of uniform convergence on
compacts as well as the corresponding Borel
$\s$-algebra $\B$. In keeping with
the literature, elements of $\B$ are called
\emph{events}.

Denote by $\mu$
the Wiener measure on $(C(\R_+),\B)$.
Recall that an event $\Lambda$ is said to hold almost surely [a.s.] if
$\mu(\Lambda)=1$.

Next we define
$U:=\{U_s\}_{s\ge 0}$ to be the \emph{Ornstein--Uhlenbeck
process} on $C(\R_+)$. The process $U$ is characterized
by the following requirements:
\begin{enumerate}
   \item It is a stationary
      infinite-dimensional diffusion with value in $C(\R_+)$;\par
   \item Its invariant measure is $\mu$.
      This implies that for any fixed $s\ge 0$, $\{U_s(t)\}_{t\ge 0}$
      is a standard linear Brownian motion.
   \item For any given $t\ge 0$, $\{U_s(t)\}_{s\ge 0}$ is a standard
      Ornstein--Uhlenbeck process on $\R$; i.e., it satisfies the stochastic
      differential equation,
      \begin{equation}\label{ou-R}
         dU_s(t) = -U_s(t)\, ds +\sqrt{2}\, dX_s\qquad{}^\forall s\ge 0,
      \end{equation}
      where $X$ is a Brownian motion.
\end{enumerate}

Following P.~Malliavin~\ycite{malliavin}, we say that
an event $\Lambda$ holds
\emph{quasi-surely} [q.s.] if
\begin{equation}\label{def:qs}
   \P\left\{ U_s\in\Lambda\text{ for all }s\ge 0\right\}=1.
\end{equation}
Because $t\mapsto U_s(t)$ is a Brownian motion,
any event $\Lambda$ that holds q.s.\@
also holds a.s.
The converse is not always true. For example,
define $\Lambda_0$ to be the collection of
all functions $f\in C(\R_+)$ that satisfy $f(1)\neq 0$
\cite{fukushima}. Evidently, $\Lambda_0$ holds a.s.\@
because with probability one Brownian motion at time
one is not at the origin.
On the other hand,  $\Lambda_0$ does not
hold q.s.\@ because $\{U_s(1)\}_{s\ge 0}$ is
point-recurrent. So the chances are $100\%$ that $U_s(1)=0$
for some $s\ge 0$.

Despite the preceding disclaimer,
a number of interesting classical
events of full Wiener measure do
hold q.s. A notable example is
a theorem of M.~Fukushima~\ycite{fukushima}. We
can state it,
somewhat informally, as follows:
\begin{equation}\label{fukushima}
   \text{The Law of the Iterated Logarithm (LIL) of
   Khintchine~\ycite{khintchine}
   holds q.s.}
\end{equation}
It might help to recall Khintchine's theorem:
For $\mu$-every $f\in C(\R_+)$,
\begin{equation}\label{lil}
   \limsup_{t\to\infty}
   \frac{f(t)}{\sqrt{2t\ln\ln t}} =1.
\end{equation}
Thus we are led to the
precise formulation of
(\ref{fukushima}): With probability one,
the continuous function $f:=U_s$ satisfies
(\ref{lil}), simultaneously for all $s\ge 0$.

For another example consider ``the other LIL'' which was
discovered by K.~L.\@ Chung~\ycite{chung}. Chung's LIL states that
for $\mu$-almost every $f\in C(\R_+)$,
\begin{equation}\label{chung}
   \liminf_{t\to\infty} \frac{\sup_{u\in[0,t]}|f(u)|}{
   \sqrt{t/\ln\ln t}} =\frac{\pi}{\sqrt 8}.
\end{equation}
Fukushima's method can be adapted to prove that
\begin{equation}\label{chung:fukushima}
   \text{Chung's LIL holds q.s.}
\end{equation}
To be more precise: With probability one,
the continuous function $f:=U_s$ satisfies
(\ref{chung}) simultaneously for all $s\ge 0$.

T.~S.~Mountford~\ycite{mountford} has derived the
quasi-sure integral test corresponding
to (\ref{fukushima}). One of the
initial aims of this article was to complement
Mountford's theorem by finding a
precise quasi-sure integral test
for~(\ref{chung:fukushima}). Before presenting
this work, let us introduce the notion
of ``relative capacity.''

For all Borel sets $G\subseteq\R_+$ and
$\Lambda\in C(\R_+)$ define
\begin{equation}
   \C_G(\Lambda) := \int_0^\infty
   \P\left\{ U_s\in\Lambda\text{ for some }s\in
   G\cap [0,\s]\right\}
   e^{-\s}\, d\s.
\end{equation}
We think of $\C_G(\Lambda)$ as the
\emph{capacity of $\Lambda$ relative to the
coordinates in $G$.}
The special case $\C_{\R_+}$ is well known
and well studied~\cite{fukushima};
$\C_{\R_+}$ is called \emph{the capacity on Wiener space}.
According to~(\ref{def:qs}), an event $\Lambda$
holds q.s.\@ iff its complement has zero $\C_{\R_+}$-capacity.

The case where
$G :=\{s\}$ is a singleton is even better studied
because of the simple fact that
$\C_{\{s\}}$ is a multiple of the Wiener measure.
Thus, $G\mapsto \C_G(\Lambda)$
interpolates from the Wiener measure ($G=\{0\}$) to
the standard capacity on Wiener space ($G=\R_+$).
This ``interpolation'' property was announced in the Abstract.

Now let $H:\R_+\to\R_+$ be decreasing and measurable, and
define
\begin{equation}
   \mathscr{L} (H) := \left\{
   f\in C(\R_+):\,
   \liminf_{t\to\infty} \left[ \sup_{u\in[0,t]}
   \;|f(u)| - H(t)\sqrt{t}\;\right] > 0
   \right\}.
\end{equation}

A decreasing measurable function $H:\R_+\to\R_+$
is called an \emph{a.s.-lower function} if
$\mathscr{L}(H)$ holds a.s.; i.e.,
$\mu$-almost every $f\in C(\R_+)$ is in $\mathscr{L}(H)$.
Likewise, $H$ is called a \emph{q.s.-lower function}
if $\mathscr{L}(H)$ holds q.s.
[The literature actually calls the function $t \mapsto H(t) \sqrt{t}$ an
a.s.[q.s]-lower function if $\mathscr{L}(H)$ holds a.s.[q.s.],
but we find our parameterization here convenient.]

To understand the utility of these definitions  better,
consider the special case that
$H(t)=\sqrt{c/\ln\ln t}$
for a fixed $c>0$ ($t\ge 0$). In this case, Chung's
LIL~(\ref{chung}) states that $\mathscr{L}(H)$
holds a.s.\@ if $c<\pi/\sqrt 8$;
its complement holds a.s.\@ if $c>\pi/\sqrt 8$.
In fact, a precise $\P$-a.s.\@ integral test is known~\cite{chung};
see Corollary~\ref{cor:chung} below.

We aim to characterize exactly when
$(\mathscr{L}(H))^\complement$ has positive $\C_G$-capacity.
Define $\K_G$ to be the
\emph{Kolmogorov $\e$-entropy} of $G$~\cite{dudley,tihomirov};
i.e., for any $\e>0$,  $k= \K_E(\e)$ is the maximal number of points
$x_1,\ldots,x_k\in E$ such that whenever $i\neq j$,
$|x_i-x_j|\ge\e$.


\begin{theorem}\label{thm:main}
   Choose and fix a decreasing
   measurable function $H:\R_+\to\R_+$, and
   a bounded Borel set $G\subset\R_+$. Then,
   $\C_G((\mathscr{L}(H))^\complement)=0$
   if and only if there exists a
   decomposition $G=\cup_{n=1}^\infty G_n$ in
   terms of closed sets $\{G_n\}_{n=1}^\infty$,
   such that
   \begin{equation} \label{eq:main_thm}
      \int_1^\infty \frac{\K_{G_n} (H^6(s))
      }{s H^2(s)}\exp\left(-\frac{\pi^2}{8 H^2(s)} \right)\, ds
      <\infty\qquad{}^\forall n\ge 1.
   \end{equation}
\end{theorem}

Theorem~\ref{thm:main} yields the following definite refinement
of~(\ref{chung}).

\begin{corollary}\label{cor:main}
   Choose and fix a decreasing
   measurable function $H:\R_+\to\R_+$. Then,
   $\mathscr{L}(H)$ holds q.s.\@ if and only if
   \begin{equation}
      \int_1^\infty \exp
      \left(-\frac{\pi^2}{8 H^2(s)} \right)\frac{ds
      }{s H^8(s)} <\infty.
   \end{equation}
\end{corollary}
Theorem~\ref{thm:main} also
contains the original almost-sure integral test
of Chung~\ycite{chung}. To prove this, simply plug
$G=\{u\}$ in Theorem~\ref{thm:main}. Then,
$\K_{\{u\}\cap J}(\e)$ is one if $u\in J$
and zero otherwise. Thus we obtain the following.

\begin{corollary}[\ocite{chung}]\label{cor:chung}
   Choose and fix a decreasing
   measurable function $H:\R_+\to\R_+$. Then
   $\mathscr{L}(H)$ holds a.s.\@ if and only if
   \begin{equation} \label{eq:chung_it}
      \int_1^\infty \exp\left(-\frac{\pi^2}{8 H^2(s)} \right)\frac{ds
      }{s H^2(s)} <\infty.
   \end{equation}
\end{corollary}

To put the preceding in perspective define
\begin{equation}
   H_\nu (t) := \frac{\pi}{\sqrt{8\left(
   \ln_+\ln_+ t + \nu\ln_+\ln_+\ln_+ t \right)}}
   \qquad{}^\forall t,\nu>0.
\end{equation}
[$1/0:=\infty$]
Then, we can deduce from Corollaries~\ref{cor:main}
and~\ref{cor:chung} that $\mathscr{L}(H_\nu)$ occurs
q.s.\@ iff $\nu>5$, whereas $\mathscr{L}(H_\nu)$ occurs
a.s.\@ iff $\nu>2$. In particular, $\mathscr{L}(H_\nu)$ occurs
a.s.\@ but not q.s.\@  if
$\nu\in[2,5)$. The following is another interesting
consequence of Theorem~\ref{thm:main}.

\begin{corollary}\label{cor:pdim}
   Let $G\subseteq[0,1]$ be a non-random Borel set.
   Then,
   \begin{equation}\begin{split}
      \Pdim G>\frac{\nu-2}{3}\ &\Longrightarrow\
         \C_G\left( \left(\mathscr{L}(H_\nu) \right)^\complement
         \right)>0,\text{ whereas}\\
      \Pdim G<\frac{\nu-2}{3}\ &\Longrightarrow\
         \C_G\left( \left(\mathscr{L}(H_\nu) \right)^\complement
         \right)=0.
   \end{split}\end{equation}
   Here, $\Pdim G$ denotes the packing dimension~\cite{Mattila}
   of the set $G$.
\end{corollary}

Throughout this paper, uninteresting constants
are denoted by $a$, $b$, $\alpha$, $A$, etc.
Their values may change from line to line.\medskip

\noindent\textbf{Acknowledgements.}\
We wish to thank Professor
Zhan Shi for generously sharing with
us the English translation of~\ocite{LS}.

\section{Brownian Sheet, and
   Capacity in Wiener Space}\label{sec:main}

We will be working with a special construction of the process
$U$. This construction is due
to D.~Williams~\cite{meyer}*{Appendix}.

Let $B:=\{B(s,t)\}_{s,t\ge 0}$ denote a two-parameter
Brownian sheet. This means that $B$ is a centered, continuous,
Gaussian process with
\begin{equation}
   \mathrm{Cov} \left( B(s,t) ~,~ B(s',t') \right) =
   \min\left( s,s'\right) \times \min\left( t,t'\right)
   \qquad{}^\forall s,s',t,t'\ge 0.
\end{equation}
The Ornstein--Uhlenbeck process
$U=\{U_s\}_{s\ge 0}$ on $C(\R_+)$ is
precisely the infinite-dimensional
process that is defined by
\begin{equation}\label{U}
   U_s(t) = \frac{B(e^s,t)}{e^{s/2}}
   \qquad{}^\forall s,t\ge 0.
\end{equation}
Indeed, one can check directly that $U$ is a
$C(\R_+)$-valued,
stationary, symmetric diffusion. And
for every $t\ge 0$, $\{U_s(t)\}_{s\ge 0}$
solves the stochastic
differential equation~(\ref{ou-R}) of the Ornstein--Uhlenbeck type.
Furthermore, the invariant
measure of $U$ is the Wiener measure.

The following well--known result
is a useful localization tool.

\begin{lemma}\label{lem:localization}
   For all bounded Borel sets $G\subseteq\R_+$
   and $\Lambda\in\B$,
   $\C_G (\Lambda)>0$ iff with positive probability
   there exists $s\in G$ such that $U_s\in\Lambda$.
\end{lemma}

\begin{remark}
   The previous lemma continues to hold even when $G$ is unbounded.
\end{remark}

\begin{proof}
   Without loss of much generality, we may---and will---assume
   that $G\subseteq[0,q]$ for some $q>0$.
   Let $p_G (\Lambda)$ denote the probability that
   there exists $s\in G$ such that $U_s\in\Lambda$.
   Evidently, $\C_G(\Lambda)\le p_G(\Lambda)$. Furthermore,
   $\C_G (\Lambda) = \int_0^q \P\{
   {}^\exists  s\in G\cap [0,\tau]:\ U_s\in\Lambda\}
   e^{-\tau} \, d\tau + e^{-q} p_G (\Lambda)$, whence
   the bounds,
   \begin{equation}\label{capP}
      e^{-q} p_G (\Lambda) \le \C_G (\Lambda) \le p_G (\Lambda).
   \end{equation}
   The lemma follows.
\end{proof}

Define
\begin{equation}\label{f*}
   f^* := \sup_{u\in[0,1]} |f(u)|\qquad{}^\forall
   f\in C(\R_+).
\end{equation}
The following is the main step in the proof of
Theorem~\ref{thm:main}. It
was announced earlier in the Abstract.

\begin{theorem}\label{main-EST}
   There exists $a>1$ such that
   for all $r\in(0,1)$
   and all Borel sets $G\subseteq[0,1]$,
   \begin{equation}
      \frac 1a \K_G (r^6) e^{-\pi^2/(8r^2)}
      \le\C_G \left\{ f^*\le r\right\}
      \le a \K_G (r^6) e^{-\pi^2/(8r^2)} .
   \end{equation}
\end{theorem}

\begin{remark}
   The constant $a$ depends on $G$ only through the fact that
   $G$ is a subset of $[0,1]$. Therefore, there exists
   $a>1$ such that simultanously for all Borel sets
   $F,G\subseteq[0,1]$,
   \begin{equation}
      \frac1a \frac{\K_F(r^6)}{\K_G(r^6)} \le
      \frac{\C_F\left\{ f^*\le r\right\}}{\C_G
      \left\{ f^*\le r\right\}} \le
      a\frac{\K_F(r^6)}{\K_G(r^6)}\qquad{}^\forall
      r\in(0,1).
   \end{equation}
\end{remark}

\begin{remark}\label{rem:cap[0,1]}
   It turns out that for any fixed $\e>0$,
   $\C_{\R_+}$ and $\C_{[0,\e]}$ are
   equivalent. To prove this, we can assume
   without loss of generality that $\e\in(0,1)$.
   [This is because $\e\mapsto\C_{[0,\e]}(\Lambda)$
   is increasing.]
   Now, on one hand, $\C_{[0,\e]}(\Lambda)\le \C_{\R_+}(\Lambda)$.
   On the other hand,
   \begin{equation}\begin{split}
      \C_{\R_+}(\Lambda) & \le \int_0^\infty
         \sum_{0\le j\le \s/\e}
         \P \left\{ {}^\exists s\in[j\e,(j+1)\e]:\
         U_s\in\Lambda \right\} e^{-\s}\, d\s\\
      & \le \P \left\{ {}^\exists s\in[0,\e]:\
         U_s\in\Lambda \right\}
         \int_0^\infty \frac{\s+1}{\e}
         e^{-\s}\, d\s,
   \end{split}\end{equation}
   by stationarity. In the notation of Lemma~\ref{lem:localization},
   the last term is $(2/\e)p_{[0,\e]}(\Lambda)\le
   (2e/\e)\C_{[0,\e]}(\Lambda)$; cf.\@ (\ref{capP}). Thus,
   \begin{equation}
      \frac{\e}{2e} \C_{\R_+}(\Lambda)\le
      \C_{[0,\e]}(\Lambda) \le \C_{\R_+}(\Lambda)
      \qquad{}^\forall\Lambda\in\B.
   \end{equation}
   This proves amply the claimed equivalence of
   $\C_{[0,\e]}$ and $\C_{\R_+}$.
\end{remark}

According to the eigenfunction
expansion of Chung~\ycite{chung},
\begin{equation}\label{eq:Chung-f}
   \mu \left\{ f^*\le r\right\} \sim \frac{4}{\pi}
   e^{-\pi^2/(8r^2)} \qquad(r\to 0).
\end{equation}
Therefore,
thanks to~(\ref{capP}), Theorem~\ref{main-EST}
is equivalent to our next result.

\begin{theorem}\label{thm:main-EST}
   Recall that $U^*_s =\sup_{t\in[0,1]}|U_s(t)|$
   [eq.\@ (\ref{f*})].
   Then, there exists a constant $a>1$ such that for
   all $r\in(0,1)$ and all Borel sets $G\subseteq[0,1]$,
   \begin{equation}
      \frac 1a\K_G (r^6) \mu \left\{ f^*\le r\right\} \le
      \P\left\{ \inf_{_{\scriptstyle
      s\in G}} U^*_s  \le r\right\}
      \le a \K_G (r^6) \mu \left\{ f^*\le r\right\}.
   \end{equation}
\end{theorem}
We will derive this particular reformulation of
Theorem~\ref{main-EST}.
The following result
plays a key role in our analysis.

\begin{proposition}[%
   Lifshits and Shi~\ycite{LS}*{Proposition 2.1}]\label{pr:SL}
   Let $\{X_t\}_{t\ge 0}$ denote planar Brownian motion.
   For every $r>0$ and $\lambda\in(0,1]$
   define
   \begin{equation}
      \mathscr{D}_\lambda^r  = \left\{ (x,y)\in\R^2:\
      |x|\le r ~,~ \left| x\sqrt{1-\lambda}
      +y\sqrt\lambda\right|\le r \right\}.
   \end{equation}
   Then there exists an $a\in(0, 1/2)$
   such that for all $r>0$ and $\lambda\in(0,1]$,
   \begin{equation}
      \P\left\{ X_t\in
      \mathscr{D}_\lambda^r \quad {}^\forall
      t\in[0,1] \right\} \le \frac1a \mu \left\{ f^*\le r\right\}
      e^{-a\lambda^{1/3}/r^2}.
   \end{equation}
\end{proposition}
%
%
\begin{lemma}\label{lem:corr}
   There exists
   a constant $a\in(0,1)$ such that for
   all $1\ge S>s>0$,
   \begin{equation}\label{eq:corr}
      \P\left\{ U^*_s \le r ~,~ U^*_S\le r\right\}
      \le \frac1a \mu \left\{ f^*\le r\right\}
      e^{- a (S-s)^{1/3}/r^2} \quad{}^\forall r\in(0,1).
   \end{equation}
\end{lemma}

\begin{proof}
   Define $\lambda=1-e^{-(S-s)}$. Then
   owing to (\ref{U}) we can write
   \begin{equation}\label{eq:decompose}
      U_S(t)   = U_s(t) \sqrt{1-\lambda} +
      \frac{B(e^S,t)-B(e^s,t)}{\sqrt{e^S-e^s}}
      \sqrt{\lambda} := U_s(t) \sqrt{1-\lambda} +
      V(t) \sqrt{\lambda}.
   \end{equation}
   By the Markov properties of the Brownian sheet,
   $X_t:=(U_s(t),V(t))$ defines a planar Brownian motion. Moreover,
   $\P\{ U^*_s \le r ~,~ U^*_S \le r \}
   =\P\{ X_t \in \mathscr{D}_\lambda^r,\ {}^\forall t\in[0,1]\}$.
   By Taylor's expansion,
   $1-e^{-x} \ge (x/2)$ ($x\in[0,1]$). Therefore,
   Proposition~\ref{pr:SL} completes the proof.
\end{proof}

\begin{proof}[Proof of Theorem~\ref{thm:main-EST}:
   Lower Bound]
   Let $k=\K_G(r^6)$, and choose maximal Kolmogorov points
   $s(1)<\cdots < s(k)$ such that
   $s(i+1)-s(i)\ge r^6$. Evidently, whenever
   $j>i$ we have
   $s(j)-s(i) \ge (j-i)r^6$. Now define
   \begin{equation}\label{OU}
      N_r = \sum_{i=1}^k \mathbf{1}_{
      \{ U^*_{s(i)} \le r \}}.
   \end{equation}
   According to Lemma~\ref{lem:corr},
   \begin{equation}\label{eq:2nd}\begin{split}
      \E\left[ N^2_r \right] &= k \mu \left\{ f^*\le r\right\}
         + 2\sum_{i=1}^{k-1}\sum_{j=i+1}^k \P\left\{
         U^*_{s(i)} \le r ~,~ U^*_{s(j)} \le r \right\}\\
      & \le k \mu \left\{ f^*\le r\right\}
         + \frac{2}{a}  \mu \left\{ f^*\le r\right\}
         \sum_{i=1}^{k-1}\sum_{j=i+1}^k
         \exp\left( - \frac{a
         (s(j)-s(i))^{1/3}}{r^2}\right)\\
      & \le k \mu \left\{ f^*\le r\right\} + \frac{2
         }{a} \mu \left\{ f^*\le r\right\}
         \sum_{i=1}^{k-1}\sum_{j=i+1}^k
         \exp\left( -a(j-i)^{1/3}\right)\\
      & \le  A k\mu \left\{ f^*\le r\right\}.
   \end{split}\end{equation}
   Note that $A$ is a positive and finite constant that does
   not depend on $r$. Also note that
   $\E[N_r]=k\mu\{ f^*\le r\}$.
   This and the Paley--Zygmund
   inequality~\cite{Khoshnevisan}*{Lemma 1.4.1, p.\@ 72} together
   reveal that
   \begin{equation}
      \P\left\{ \inf_{s\in G} U^*_s\le r \right\} \ge
      \P\left\{ N_r >0 \right\}
      \ge \frac{\left( \E[N_r] \right)^2}{\E\left[ N^2_r\right]}
      \ge \frac{k}{A}\mu \left\{ f^*\le r\right\}.
   \end{equation}
   The definition of $k$
   implies the lower bound in Theorem~\ref{thm:main-EST}.
\end{proof}

Before proving the upper bound
of Theorem~\ref{thm:main-EST} in
complete generality,
we first derive the following weak form:

\begin{proposition}\label{pr:main-est}
    There exists a finite constant $a>1$
    such that for all $r\in(0,1)$,
    $\P\{ \inf_{s\in[0,r^6]} U^*_s \le r \}
    \le a \mu \left\{ f^*\le r\right\}$.
\end{proposition}

\begin{proof}
   Recall~(\ref{OU}), and define
   \begin{equation}
      L  (s;r) = \int_0^s \mathbf{1}_{\{ U^*_\nu\le r\}}\, d\nu
      \qquad{}^\forall s,r>0.
   \end{equation}
   Let $\F:=\{\F_s\}_{s\ge 0}$ denote
   the augmented filtration generated by
   the infinite-dimensional process $\{U_s\}_{s\ge 0}$.
   The latter process is Markov with respect to $\F$. Moreover,
   \begin{equation}
      \E\left[ \left. L  (2r^6; r+r^3 )\, \right|\,
      \F_s \right]  \ge \int_s^{2r^6}
      \P\left\{ \left.
      U^*_\nu\le  r+r^3 \, \right|\, \F_s\right\}\, d\nu
      \cdot \mathbf{1}_{\{ U^*_s \le r\}}.
   \end{equation}
   As in (\ref{eq:decompose}),
   if $\nu>s$ are fixed, then we can write
   \begin{equation}\begin{split}
      U_\nu(t) & = U_s(t) e^{-(\nu-s)/2} +
            \frac{B(e^\nu,t)-B(e^s,t)}{\sqrt{e^\nu-e^s}}
            \sqrt{1-e^{-(\nu-s)}}\\
         & := U_s(t) e^{-(\nu-s)/2} + V(t) \sqrt{1-e^{-(\nu-s)}}.
   \end{split}\end{equation}

   We emphasize, once again, that $(U_s,V)$ is a
   planar Brownian motion. In addition, $V$ is
   independent of $\F_s$, and
   $U_\nu^*  \le U_s^* + V^* \sqrt{1-\exp\{-(\nu-s)\}}.$
   Consequently,  as long as
   $0\le s\le r^6$ and $s<\nu<2r^6$,
   \begin{equation}\label{trick}
      U_\nu^*  \le U_s^* + \frac{ r^3 }{\sqrt 2} V^*.
   \end{equation}
   [We have used the inequality $1-e^{-z} \le z/2$ valid for
   all $z\in(0,1)$.]
   Therefore, for all $0\le s\le r^6$,
   \begin{equation}\begin{split}
      M(s) &= \E\left[ \left. L  (2r^6; r+r^3 )\, \right|\,
         \F_s \right] \\
      & \ge \int_s^{2r^6} \P\left\{
         V^* \le \sqrt 2 \right\}\, d\nu
         \cdot \mathbf{1}_{\{ U^*_s\le r \}}\\
      &  =\mu \left\{ f^*\le \sqrt 2 \right\}
         \left(2r^6 -s\right)
         \cdot \mathbf{1}_{\{ U^*_s \le r\}}\\
      & \ge \mu \left\{ f^*\le \sqrt 2 \right\} r^6
         \cdot \mathbf{1}_{\{ U^*_s \le r\}}.
   \end{split}\end{equation}
   Because $\{M(s)\}_{s\ge 0}$ is a martingale,
   we can apply Doob's
   maximal inequality to obtain the following:
   \begin{equation}\begin{split}
      \P\left\{ \inf_{s\in[0,r^6]} U^*_s \le r\right\}
         & \le \P\left\{ \sup_{s\in[0,r^6]} M(s) \ge
         \mu \left\{ f^*\le \sqrt 2\right\} r^6  \right\}\\
      & \le \frac{\E\left[ L  (2r^6; r+r^3 )\right]}{
         \mu \left\{ f^*\le \sqrt 2\right\} r^6}
         = \frac{2\mu\left\{ f^*\le r+r^3 \right\}}{
         \mu \left\{ f^*\le \sqrt 2\right\}}.
   \end{split}\end{equation}
   Thanks to (\ref{eq:Chung-f}),
   \begin{equation}
      \frac{\mu\left\{ f^*\le r+r^3\right\}}{
      \mu\left\{ f^*\le r\right\}}  \sim
      \exp\left( - \frac{\pi^2}{8} \left[
      \frac{1}{( r+r^3 )^2}- \frac{1}{r^2}
      \right]\right) \to e^{\pi^2/4}.
      \qquad( r\to  0).
   \end{equation}
   Thus, the left-hand side is bounded
   ($r\in(0,1)$), and the proposition follows.
\end{proof}

\begin{proof}[Proof of Theorem~\ref{thm:main-EST}: Upper Bound]
   Define $n=n(r)$ to be $\lfloor r^{-6} \rfloor$, and define
   $I(j;n)$ to be the interval $[j/n,(j+1)/n)$
   $(j=0,\ldots,n)$. Then, by stationarity and
   Proposition~\ref{pr:main-est},
   \begin{equation}\label{eq:stronger}
      \P\left\{ \inf_{s\in G} U^*_s  \le r\right\}
      \le \sum_{\scriptstyle 0\le j\le n:\
      \atop\scriptstyle I(j;n)\cap G \neq \varnothing }
      \P\left\{ \inf_{s\in I(j;n)} U^*_s \le r \right\} \le a
      \mu\left\{ f^*\le r\right\} \M_n(G),
   \end{equation}
   where $ \M_n(G)=\#\{
   0\le j\le n:\ I(j;n)\cap G\neq \varnothing \}$ defines the
   \emph{Minkowski content} of $G$. In the companion to this
   paper~\ycite{KLM2}*{Proposition 2.7}  we proved that
   $\M_n(G)\le 3\K_G(1/n)$.
   By monotonicity, the latter is at most $3\K_G(r^6)$,
   whence the theorem.
\end{proof}

\section{Proof of Theorem~\ref{thm:main}
   and Corollaries~\ref{cor:main} and~\ref{cor:pdim}}

We begin with some preliminary discussions. Define
\begin{equation}\label{eq:s}
   \psi_H(G) := \int_1^\infty \frac{\K_G(H^6(s))}{s
   H^2(s)}\exp\left( -
   \frac{\pi^2}{8H^2(s)} \right)\, ds,\qquad
   \s(r) := \mu\left\{ f^*\le r\right\}.
\end{equation}

Following Erd\H{o}s~\ycite{erdos}, define
\begin{equation}
   \ee_n = e^{n \left/ \ln_+ n\right.},
   \quad H_n = H(\ee_n)\qquad{}^\forall n\ge 1.
\end{equation}

The ``critical'' function in (\ref{eq:chung_it})
is $H^2(t)=\pi^2/(8\ln_+\ln_+  t)$. This,
the fact that $\pi/\sqrt 8\in(1,2)$, and
a familiar argument~\cite{erdos}*{equations 1.2 and 3.4},
together allow us to assume
without loss of generality that
\begin{equation}\label{eq:wlog}
   \frac{1}{\sqrt{\ln_+ n}}
   \le H_n \le \frac{2}{\sqrt{\ln_+ n}}
   \qquad{}^\forall n\ge 1.
\end{equation}
From this we can conclude the existence of
a constant $a>1$
such that
\begin{equation}\label{eq:key-ee}
   \frac 1a H_n^2 \ee_{n+1} \le
   \ee_{n+1} - \ee_n \le a
   H_{n+1}^2 \ee_n\qquad{}^\forall n\ge 1.
\end{equation}
According to our companion work~\ycite{KLM2}*{eq. 2.8}, for
all $r>0$ sufficiently small,
\begin{equation}\label{eq:KGKG}
   \K_G(\e)\le 6\K_G(2\e).
\end{equation}

Because $\ee_{n+1}\sim\ee_n$ as $n\to\infty$,
(\ref{eq:Chung-f}), (\ref{eq:key-ee}), and~(\ref{eq:KGKG})
together imply that
\begin{equation}\label{eq:sumint}
    \sum_{n=1}^\infty \K_G \left( H^6_n \right)
    \s\left(H_n \right)<\infty\quad
    \text{iff}\quad \psi_H(G) <\infty.
\end{equation}

The following is the key step toward proving Theorem~\ref{thm:main}.

\begin{proposition}\label{pr:chung1}
   Let $H:\R_+\to\R_+$ be decreasing and measurable.
   Then for all non-random Borel sets $G\subseteq[0,1]$,
   \begin{equation}\begin{split}
      \liminf_{_{\scriptstyle t\to\infty}}
      \left( \inf_{_{\scriptstyle s\in G}}
      \sup_{u\in[0,t]} |U_s(u)| - H(t)\sqrt{t}
      \right) =
      \begin{cases}
         +\infty, & \text {if } \psi_H(G)<\infty,\\
         -\infty, & \text{if } \psi_H(G)= \infty.
      \end{cases}
   \end{split}\end{equation}
\end{proposition}

First we assume this proposition
and derive Theorem~\ref{thm:main}. Then,
we will tidy things up by proving
the technical Proposition~\ref{pr:chung1}.

Let us recall~(\ref{eq:s}).

\begin{definition}
   We say that $\Psi_H(G)<\infty$ if we can decompose
   $G$ as $G=\cup_{n=1}^\infty G_n$---where
   $G_1,G_2,\ldots$ are closed---such that for all
   $n\ge 1$, $\psi_H(G_n)<\infty$.
   Else, we say that $\Psi_H(G)=\infty$.
\end{definition}

Let us first rephrase Theorem~\ref{thm:main} in
the following convenient, and equivalent, form.

\begin{proposition}\label{pr:chung2}
   Let $H:\R_+\to\R_+$ be decreasing and measurable
   and $G\subseteq[0,1]$ be non-random and Borel.
   If $\Psi_H(G)<\infty$ , then
   \begin{equation}\begin{split}
      \inf_{s\in G} \liminf_{t\to\infty}
      \left( \sup_{u\in[0,t]} |U_s(u)| - H(t)\sqrt{t}
      \right) = \infty\qquad\P\text{-a.s.}
   \end{split}\end{equation}
   Else, the left-hand side is $\P$-a.s.\@ equal to $-\infty$.
\end{proposition}

\begin{proof}[Proof of Theorem~\ref{thm:main}
   in the form of Proposition~\ref{pr:chung2}]
   First suppose $\Psi_H(G)$ is finite.  We can write
   $G=\cup_{n=1}^\infty G_n$, where the $G_n$'s are
   closed and $\psi_H(G_n)<\infty$ for all $n\ge 1$.
   Then, according
   to Proposition~\ref{pr:chung1},
   \begin{equation}\begin{split}
      &\inf_{s\in G_n}\liminf_{t\to\infty}
         \left[ \sup_{u\in[0,t]} |U_s(u)| - H(t)\sqrt{t} \right]\\
      & \ge \liminf_{t\to\infty}\inf_{s\in G_n}
         \left[ \sup_{u\in[0,t]} |U_s(u)| - H(t)\sqrt{t} \right]
         =\infty.
   \end{split}\end{equation}
   This proves that
   $\inf_{s\in G}\liminf_{t\to\infty}
   ( \sup_{u\in[0,t]} |U_s(u)| - H(t)\sqrt{t} )=\infty$ a.s.\@ [$\P$].

   For the converse portion
   suppose $\Psi_H(G)=\infty$,
   and choose arbitrary non-random
   closed sets $\{G_n\}_{n=1}^\infty$ such that $\cup_{n=1}^\infty G_n=G$.
   By definition, $\psi_H(G_n)=\infty$ for some $n\ge 1$.
   Define for all $T \ge 1$,
   \begin{equation}
     \mathscr{S}_T := \left\{ s\in [0,1]:\ \inf_{t \geq T}
      \frac{\sup_{u\in[0,t]}
      \left| U_s(u) \right|}{ H(t)\sqrt{t}} \leq 1
      \right\}.
   \end{equation}
   Evidently, $\mathscr{S}_T$ is a random set for each $T \ge 0$.
   Moreover, the continuity of the Brownian sheet implies
   that with probability one, $\mathscr{S}_T$ is closed
   for all $T$; hence, so is $\mathscr{S}_T \cap G_n$.   Because
   $\psi_H(G_n)=\infty$, Proposition~\ref{pr:chung1} implies that almost
   surely,
   $\mathscr{S}_T \cap G_n \neq \varnothing$.
   Since $\{ \mathscr{S}_T \cap G_n \}_{T = 1}^\infty $ is a decreasing
   sequence of non-void compact sets, they have non-void intersection.
   That is, $(\cap_{T = 1}^\infty \mathscr{S}_T) \cap G_n
   \neq\varnothing$ a.s.\@ [$\P$].
   Replace $H$ by $H-H^3$ to complete the proof of
   Proposition~\ref{pr:chung2}.
\end{proof}

Now we derive Proposition~\ref{pr:chung1}.
This completes our proof of
Theorem~\ref{thm:main}. Our proof is divided naturally
into two halves.

\begin{proof}[Proof of Proposition~\ref{pr:chung1}: First Half]
   Throughout this portion of the proof,
   we assume that $\psi_H(G)<\infty$.

   Because $\ee_{n+1}\sim\ee_n$ as $n\to\infty$,
   Theorem~\ref{thm:main-EST} and Brownian scaling together
   imply that
   \begin{equation}\begin{split}
      &\P\left\{ \inf_{_{\scriptstyle s\in G}} \sup_{u\in[0,\ee_{n-1}]}
         |U_s(u)| \le H_n \sqrt{\ee_n} \right\}  =
         \P\left\{ \inf_{_{\scriptstyle s\in G}}
         U^*_s \le H_n \sqrt{\ee_n/\ee_{n-1}} \right\}\\
      &\le a \K_G\left(
         H_n^6 \left[ \frac{\ee_n}{\ee_{n-1}} \right]^3 \right)
         \s \left(
         H_n \sqrt{ \frac{\ee_n}{\ee_{n-1}} } \right).
   \end{split}\end{equation}
   According to (\ref{eq:KGKG}),
   $\K_G(\cdots)\le 6 \K_G(H_n^6)$ 
   for all $n$ large. This and (\ref{eq:key-ee})
   together imply that for all $n$ large,
   \begin{equation}\begin{split}
      &\P\left\{ \inf_{_{\scriptstyle s\in G}} \sup_{u\in[0,\ee_{n-1}]}
         |U_s(u)| \le H_n \sqrt{\ee_n} \right\}  \\
      &\le a \K_G \left( H^6_n \right)
         \s\left(
         H_n \sqrt{ 1 + A H_{n+1}^2 }\;\right)\\
      &\le a \K_G \left( H^6_n \right)
         \s\left( H_n \left[ 1 + A H_n^2
         \right] \right).
   \end{split}\end{equation}
   In accord with (\ref{eq:Chung-f}), for any fixed $c\in\R$,
   \begin{equation}\label{eq:sscale}
      \s \left( r+c r^3  \right) =O(\s(r) )\qquad( r\to  0).
   \end{equation}
   Thus, for all $n\ge 1$,
   \begin{equation}\label{eq:chung1-ub}
      \P\left\{ \inf_{_{\scriptstyle s\in G}} \sup_{u\in[0,\ee_{n-1}]}
      |U_s(u)| \le H_n \sqrt{\ee_n} \right\}
      \le a \K_G \left( H^6_n \right)
      \s\left( H_n \right).
   \end{equation}
   Because we are assuming that $\psi_H(G)$ is finite,
   (\ref{eq:sumint}) and the Borel--Cantelli lemma together
   imply that almost surely,
   $\inf_{_{\scriptstyle s\in G}} \sup_{u\in[0,\ee_{n-1}]}
   |U_s(u)| > H_n \sqrt{\ee_n}$
   for all but a finite number of $n$'s.
   It follows from this and a standard monotonicity argument that
   \begin{equation}
      \psi_H(G)<\infty\
      \Longrightarrow\
      \liminf_{_{\scriptstyle t\to\infty}} \left[
      \inf_{_{\scriptstyle s\in G}} \sup_{u\in[0,t]}
      |U_s(u)| - H(t) \sqrt{t}\right]>0\ \text{ a.s.\@ } [\P].
   \end{equation}
   But if $\psi_H(G)$ were finite then
   $\psi_{H+H^3}(G)$ is also finite; compare
   (\ref{eq:KGKG}) and~(\ref{eq:sscale}). Thanks to~(\ref{eq:wlog}),
   $\lim_{t\to\infty}H^3(t)\sqrt{t}=\infty$. Therefore,
   the $\liminf$ of the preceding display is infinity.
   This concludes the first half of our proof of
   Proposition~\ref{pr:chung1}.
\end{proof}

In order to prove the second half of Proposition~\ref{pr:chung1}
we assume that $\psi_H(G)=\infty$, recall (\ref{eq:s}),
and define
\begin{equation}\begin{split}
   L_n &:= \left\{
      \inf_{_{\scriptstyle s\in G}} \sup_{u\in[0,\ee_n]}
      \left| U_s(u) \right| \le H_n\sqrt{\ee_n}
      \right\},\\
   f(z)  &:= \K_G \left(z^6\right) \s(z).
\end{split}\end{equation}

\begin{lemma}\label{lem:PLL}
   Define for all $j\ge i$,
   $\lambda_{i,j} := \ee_j/(\ee_j-\ee_i)$ and
   $\delta_{i,j}:= H_j\sqrt{\lambda_{i,j}}
   +H_i\sqrt{\lambda_{i,j}-1}$. Then, there exists $a>1$ such that
   for all $j\ge i$,
   $\P(L_j\,|\,L_i)
   \le a\K_G \left(
   \delta_{i,j}^6\right) \s \left( \delta_{i,j} \right)$.
\end{lemma}

\begin{proof}
   Evidently, $\P(L_j\,|\,L_i)$ is at most
   \begin{equation}\begin{split}
      &\P\left\{ \left.
         \inf_{_{\scriptstyle s\in G}} \sup_{u\in[\ee_i,\ee_j]}
         \left| U_s(u) \right| \le H_j\sqrt{\ee_j}
         \ \right|\ L_i \right\}\\
      &= \P\left\{ \left.
         \inf_{_{\scriptstyle s\in G}} \sup_{u\in[\ee_i,\ee_j]}
         \left| U_s(u) - U_s(\ee_i) + U_s(\ee_i)
         \right| \le H_j\sqrt{\ee_j}
         \ \right|\ L_i \right\}\\
      &\le \P\left\{
         \inf_{_{\scriptstyle s\in G}} \sup_{u\in[\ee_i,\ee_j]}
         \left| U_s(u) - U_s(\ee_i)
         \right| \le H_j\sqrt{\ee_j} + H_i\sqrt{\ee_i} \right\}.
   \end{split}\end{equation}
   We have appealed to the Markov properties of the Brownian sheet
   in the last line. Because $u\mapsto U_\bullet(u)$ is
   a $C(\R_+)$-valued Brownian motion,
   \begin{equation}\begin{split}
      \P(L_j\,|\, L_i) &\le \P\left\{
         \inf_{_{\scriptstyle s\in G}} \sup_{u\in[0,\ee_j-\ee_i]}
         \left| U_s(u) \right| \le H_j\sqrt{\ee_j} + H_i\sqrt{\ee_i}
         \right\}\\
      & = \P\left\{
         \inf_{s\in G} U^*_s \le \delta_{i,j}\right\}.
   \end{split}\end{equation}
   Theorem~\ref{thm:main-EST} completes the proof.
\end{proof}

Our forthcoming estimates of $\P(L_j\,|\,L_i)$ rely
on the following elementary bound; see, for example,
our earlier work~\ycite{KLM1}*{eq.\@ 8.30}:
Uniformly for all integers $j>i$,
\begin{equation}\label{eq:ees}
   \ee_j-\ee_i \ge
   \ee_i\left(\frac{j-i}{\ln i}\right)(1+o(1))
   \qquad(i\to\infty).
\end{equation}

\begin{lemma}\label{c1}
   There exist $i_0\ge 1$ and
   a finite $a>1$
   such that for all $i\ge i_0$ and
   $j\ge i+\ln^{19}(j)$,
   \begin{equation}\label{eq:c1}
      \P(L_j\,|\, L_i) \le a \P(L_j).
   \end{equation}
\end{lemma}

\begin{proof}
   Thanks to~(\ref{eq:wlog}) and~(\ref{eq:ees}),
   the following holds uniformly over all $j>i+\ln^{19}(j)$:
   $(\ee_j/\ee_i) \ge (1+o(1)) H_j^{-36}$ $(i\to\infty)$.
   Thus, uniformly over all $j>i+\ln^{19}(j)$,
   \begin{equation}\begin{split}
      \sqrt{\lambda_{i,j}}
         = \frac{1}{\sqrt{1-(\ee_i/\ee_j)}}
         &\le \frac{1}{\sqrt{1-(1+o(1))H_j^{36}}}
         = 1+ O\left(H_j^3\right),\\
       H_i \sqrt{\lambda_{i,j}-1}
         &= O\left( H_j^3\right)\qquad(i\to\infty).
   \end{split}\end{equation}
   Lemma~\ref{lem:PLL} guarantees then that
   uniformly over all $j>i+\ln^{19}(j)$,
   $\delta_{i,j}\le H_j+O(H_j^3)$,
   and the big-$O$ and little-$o$ terms do not depend on
   the $j$'s in question.
   The lemma follows from this,
   equations~(\ref{eq:KGKG}) and~(\ref{eq:sscale}), and
   Theorem~\ref{thm:main-EST}.
\end{proof}

\begin{lemma}\label{c2}
   There exist $i_1\ge 1$ and $a\in(0,1)$
   such that for all $i\ge i_1$ and
   $j\in [ i+\ln (i),i+\ln^{19}(j))$,
   $\P(L_j\,|\, L_i) \le (a j^a)^{-1}.$
\end{lemma}

\begin{proof}
   Equations~(\ref{eq:ees}) and~(\ref{eq:wlog})
   together imply that uniformly
   for all $j\ge i+\ln(i)$,
   $(\ee_i/\ee_j) \le \frac12+o(1)$ $(i\to\infty)$.
   This is equivalent to the existence of a constant
   $A_{\ref{eq:AA}}$ such that for all $(i,j)$ in
   the range of the lemma,
   \begin{equation}\label{eq:AA}
      \sqrt{\lambda_{i,j}} \vee
      \sqrt{\lambda_{i,j}-1} \le a.
   \end{equation}
   Thanks to~(\ref{eq:wlog}), we can enlarge the last constant
   $a$, if necessary, to ensure that for all $(i,j)$ in the range
   of this lemma,
   $H_i\le a H_j$. Therefore,
   Lemma~\ref{lem:PLL} then implies that
   $\delta_{i,j}=O(H_j)$, and
   the big-$O$ term does not
   depend on the range of $j$'s in question.
   Because $G\subseteq[0,1]$,
   \begin{equation}\label{eq:KGKGKG}
     \K_G(\e)\le \K_{[0,1]}(\e)\sim 1/\e
     \qquad(\e\to 0).
  \end{equation}
   Thus, Lemma~\ref{lem:PLL} ensures that
   $\P(L_j\,|\,L_i)\le a
   \delta_{i,j}^{-6}\s(\delta_{i,j}).$
   Near the origin,
   the function $\delta\mapsto \delta^{-6}\s(\delta)$ is increasing.
   Because we have proved that over the range
   of $(i,j)$ of this lemma $\delta_{i,j}=O(H_j)$,
   equation (\ref{eq:Chung-f}) asserts the
   existence of a universal
   $\alpha>1$ such that
   $\P(L_j\,|\, L_i)$ is at most
   $\alpha H_j^{-6} \exp( - \alpha^{-1} H_j^{-2}).$
   Equation~(\ref{eq:wlog}) then completes our proof.
\end{proof}

\begin{lemma}\label{c3}
   There exist $i_2\ge 1$ and $a>1$
   such that for all $i\ge i_2$ and
   $j\in(i, i+\ln i)$,
   $\P(L_j\,|\, L_i) \le ae^{- (j-i)/a}.$
\end{lemma}

\begin{proof}
   By~(\ref{eq:ees}), $(\ee_i/\ee_j) \le 1 - (1+o(1))(j-i) \ln^{-1} (i)$
   $(i\to\infty)$,
   where the little-$o$ term does not depend on
   $j\in(i,i+\ln i)$. Similarly,
   $(\ee_j/\ee_i)\ge 1+(1+o(1))(j-i)\ln^{-1}(i)$.
   Thus, as $i\to\infty$,
   \begin{equation}\begin{split}
      \sqrt{\lambda_{i,j}} &=
         \frac{1}{\sqrt{1-( \ee_i/\ee_j)}}
         \le (1+o(1))\sqrt{\frac{\ln i}{j-i}}
         \le \frac{2+o(1)}{H_j\sqrt{j-i}},\\
      \sqrt{\lambda_{i,j}-1} &=
         \frac{1}{\sqrt{( \ee_j/\ee_i)-1}}
         \le (1+o(1))\sqrt{\frac{\ln i}{j-i}}
         \le \frac{2+o(1)}{H_j\sqrt{j-i}},
   \end{split}\end{equation}
   by (\ref{eq:wlog}). Once again,
   the little-$o$ terms are all independent of
   $j\in (i,i+\ln i)$. Because $H_i=O(H_j)$
   uniformly for all $(i,j)$ in the range considered
   here,  Lemma~\ref{lem:PLL} implies that
   uniformly for all $j\in (i,i+\ln i)$,
   $\delta_{i,j}=O(1/\sqrt{j-i})$.
   Equation~(\ref{eq:KGKGKG}) bounds the first term
   on the right-hand side; (\ref{eq:Chung-f}) bounds
   the second. This and~(\ref{eq:wlog})
   together prove the existence of a constant
   $\alpha>1$ such that for all $i\ge i_2$
   and all $j\in(i,i+\ln i)$,
   $\P(L_j\,|\, L_i) \le \alpha (j-i)^{3}
   \exp\{ - (j-i)/\alpha\}.$
   The lemma follows.
\end{proof}

\begin{proof}[Proof of Proposition~\ref{pr:chung1}: Second Half]
   According to Theorem~\ref{thm:main-EST},
   for all $n$ large enough,
   $\P(L_n)\ge af(H_n)$.
   Because $\psi_H(G)=\infty$, the latter estimate
   and~(\ref{eq:sumint}) together imply that
   \begin{equation}
      \sum_{i=1}^\infty \P(L_i) =\infty.
   \end{equation}
   Thus, our derivation is complete once we
   demonstrate the following:
   \begin{equation}\label{eq:EL2}
      \liminf_{n\to\infty} \frac{ \sum_{i=1}^{n-1}
      \sum_{j=i}^n \P(L_i\cap L_j)}{
      \left( \sum_{i=1}^n \P(L_i)  \right)^2}<\infty.
   \end{equation}
   See Chung and Erd\H{o}s~\ycite{ce}.
   In fact, the preceding display holds with a $\limsup$ in place
   of the $\liminf$. This fact follows from combining,
   using standard arguments, Lemmas
   \ref{c1} through~\ref{c3}.

   Indeed, let $I:=\max(3,i_1,i_2,i_3)$ and $s_n:=\sum_{i=1}^n\P(L_i)$.
   Lemma~\ref{c1} ensures that
   \begin{equation}\label{C1}\begin{split}
      \mathop{\sum_{i=I}^{n-1}\sum_{j=i}^n}\limits_{j>
      i+\ln^{19}(j)}
      \P(L_j\cap L_i) =O\left( s_n^2 \right).
   \end{split}\end{equation}

   By Lemma~\ref{c2},
   \begin{equation}\label{C2}\begin{split}
      &\mathop{\sum_{i=I}^{n-1}\sum_{j=i}^n}\limits_{
         j\in\left(i+\ln(i),i+\ln^{19}(j)\right]}
         \P(L_j\cap L_i) \le \frac 1a
         \mathop{\sum_{i=I}^{n-1}\sum_{j=i}^n}\limits_{
         j\in\left(i+\ln(i),i+\ln^{19}(j)\right]}
         j^{-a}\P(L_i)\\
      &\qquad= \sum_{i=I}^n
         O\left( \frac{\ln^{19}(i)}{i^a}
         \right) \P(L_i)
         = O\left( s_n \right).
   \end{split}\end{equation}
   The big-$O$ terms do not depend on the variables
   $(j,n)$.

   Finally, Lemma~\ref{c3} implies that
   \begin{equation}\label{C3}\begin{split}
      \mathop{\sum_{i=I}^{n-1}\sum_{j=i}^n}\limits_{j\in
         (i,i+\ln i]}
         \P(L_j\cap L_i) &\le a
         \sum_{i=1}^n \sum_{j=i}^\infty
         \P(L_i)e^{(j-i)/a} =O\left( s_n \right).
   \end{split}\end{equation}
   We have already seen that $s_n\to\infty$.
   Thus, \eqref{C1}--\eqref{C3} imply~\eqref{eq:EL2},
   and hence the theorem. More precisely, we have
   proved so far that
   \begin{equation}
      \psi_H(G)=\infty\
      \Longrightarrow\
      \liminf_{_{\scriptstyle t\to\infty}} \left[
      \inf_{_{\scriptstyle s\in G}} \sup_{u\in[0,t]}
      |U_s(u)| - H(t) \sqrt{t}\right]<0\ \text{ a.s.\@ }[\P].
   \end{equation}
   Replace $H$ by $H+H^3$ to deduce that the preceding
   $\liminf$ is in fact $-\infty$. This completes
   our proof of Proposition~\ref{pr:chung1}.
\end{proof}

We conclude this section by proving the remaining Corollaries
\ref{cor:main} and~\ref{cor:pdim}.

\begin{proof}[Proof of Corollary~\ref{cor:main}]
   By definition, $\mathscr{L}(H)$ holds q.s.\@
   iff $\C_{\R_+}((\mathscr{L}(H))^\complement)=0$.
   Thanks to Theorem~\ref{thm:main}, this condition
   is equivalent to the existence of a non-random
   ``closed-denumerable'' decomposition
   $\R_+=\cup_{n=1}^\infty G_n$ such that for all $n\ge 1$,
   $\psi_H(G_n)<\infty$. But
   one of the $G_n$'s must contain a closed interval
   that has positive length. Therefore, by the translation-invariance
   of $G\mapsto \K_G(r)$, there exists
   $\e\in(0,1)$ such that $\psi_H([0,\e])<\infty$.

   Conversely, if
   $\psi_H([0,\e])$ is finite, then we can define
   $G_n$ to be $[(n-1)\e,n\e]$ ($n\ge 1$) to find that
   $\psi_H(G_n)=\psi_H([0,\e])<\infty$. Theorem~\ref{thm:main}
   then proves that $\C_{\R_+}((\mathscr{L}(H))^\complement)=0$
   iff there exists $\e>0$ such that $\psi_H([0,\e])<\infty$.
   Because $\K_{[0,\e]}(r)\sim \e/r$ ($r\to 0$), the
   corollary follows.
\end{proof}

\begin{proof}[Proof of Corollary~\ref{cor:pdim}]
   We can change variables to deduce that
   $\psi_{H_\nu}(G)$ is finite iff
   $\int_1^\infty \K_G(1/s) s^{-1-(\nu/3)}\, ds$ converges.
   This and Proposition 2.8
   of our companion work~\ycite{KLM2} together imply
   that
   \begin{equation}
      \inf\{ \nu>0:\ \psi_{H_\nu}(G)<\infty\}
      = 2+3\Kdim G,
   \end{equation}
   where $\Kdim$ denotes the (upper) Minkowski dimension~\cite{Mattila}.
   By regularization~\cite{Mattila}*{p.\@ 81},
   \begin{equation}
      \inf\{ \nu>0:\ \Psi_{H_\nu}(G)<\infty\}
      = 2+3\Pdim G.
   \end{equation}
   Theorem~\ref{thm:main} now implies Corollary~\ref{cor:pdim}.
\end{proof}


\begin{bibdiv}
\begin{biblist}

\bib{chung}{article}{
   author = {Chung, Kai Lai},
    title = {On the maximum partial sums of sequences of independent random
            variables},
     date = {1948},
  journal = {Trans.\@ Amer.\@ Math.\@ Soc.},
   volume = {64},
    pages = {205--233},
}

\bib{ce}{article}{
    author = {Chung, K. L.},
    author = {Erd\H{o}s, P.},
     TITLE = {On the application of the {B}orel-{C}antelli lemma},
   JOURNAL = {Trans.~Amer.~Math.~Soc.},
    VOLUME = {72},
      YEAR = {1952},
     PAGES = {179--186},
}

\bib{dudley}{incollection}{
   author =   {Dudley, R.~M.},
    title =   {A Course in Empirical Processes},
booktitle =   {\'Ecole d'\'et\'e de
              probabilit\'es Saint-Flour, XII--1982},
   pages  =   {1--142},
     date =   {1973},
publisher =   {Springer},
  address =   {Berlin},
}

\bib{erdos}{article}{
   author =    {Erd\H{o}s, Paul},
    title =    {On the law of the iterated logarithm},
  journal =    {Ann.\@ Math.},
     year =    {1942},
   volume =    {43(2)},
    pages =    {419\ndash 436},
}

\bib{fukushima}{article}{
   author = {Fukushima,~Masatoshi},
    title = {Basic properties of Brownian motion
            and a capacity on the Wiener space},
  journal = {J.\@ Math.\@ Soc.\@ Japan},
   volume = {36(1)},
   pages =  {161--176},
    year =  {1984},
}

\bib{khintchine}{book}{
    AUTHOR = {Khintchine, A.~Ya.},
     TITLE = {Asymptotische {G}esetz der {W}ahrscheinlichkeitsrechnung},
 PUBLISHER = {Springer},
   ADDRESS = {Berlin},
      YEAR = {1933},
}

\bib{Khoshnevisan}{book}{
   author =   {Khoshnevisan, Davar},
    title =   {Multiparameter Processes: An Introduction
              to Random Fields},
     year =   {2002},
publisher =   {Springer},
  address =   {New York},
}

\bib{KLM1}{article}{
   author =   {Khoshnevisan, Davar},
   author =   {Levin, David A.},
   author =   {M\'endez--Hern\'andez, Pedro J.},
    title =   {On dynamical Gaussian random walks},
     year =   {2003},
  journal =   {Ann.~Probab.},
   status =   {To appear},
}

\bib{KLM2}{article}{
   author =   {Khoshnevisan, Davar},
   author =   {Levin, David A.},
   author =   {M\'endez--Hern\'andez, Pedro J.},
    title =   {Exceptional times and invariance for dynamical random
              walks},
     year =   {2004},
   status =   {Preprint},
}


\bib{LS}{article}{
    author={Lifshits, M. A.},
    author={Shi, Z.},
     title={Lower functions of an empirical process and of a Brownian sheet},
  language={Russian, with Russian summary},
   journal={Teor. Veroyatnost. i Primenen.},
    volume={48(2)},
      date={2003},
    number={2},
     pages={321\ndash 339},
}

\bib{malliavin}{article}{
  author =   {Malliavin, Paul},
  title =    {R\'egularit\'e de lois conditionnelles et calcul des
                  variations stochastiques},
  journal =  {C.R. Acad. Sci.Paris, S\'er. A-B},
  year =     {1979},
  volume =   {289},
  number =   {5},
}

\bib{Mattila}{book}{
   author = {Mattila, Pertti},
    title = {Geometry of Sets and Measures in Euclidean Spaces:
            Fractals and Rectifiability},
publisher = {Cambridge University Press},
  address = {Cambridge},
     year = {1995},
}

\bib{meyer}{article}{
  author =   {Meyer, P.-A.},
  title =    {Note sur les processus d'Ornstein--Uhlenbeck
                  (Appendice: Un resultat de D. Williams)},
  booktitle =    {S\'em. de Probab. XVI},
  pages =    {95\ndash 133},
  publisher =    {Springer},
  year =     {1982},
  volume =   {920},
  series =   {Lec. Notes in Math.},
}

\bib{mountford}{article}{
   author =   {Mountford, T.~S.},
    title =    {Quasi-everywhere upper functions},
booktitle =    {S\'em.\@ de Probab.\@ XXVI},
    pages =    {95\ndash 106},
publisher =    {Springer},
     year =    {1992},
   volume =    {1526},
   series =    {Lect.\@ Notes in Math.},
}


\bib{tihomirov}{article}{
   author =   {Tihomirov, V.~M.},
    title =   {The works of A.\@ N.\@ Kolmogorov on $\e$-entropy
              of function classes and superpositions of functions},
  journal =   {Uspehi Mat.\@ Nauk},
     year =   {1963},
   volume =   {18(5 (113))},
    pages =   {55\ndash 92},
}

\end{biblist}
\end{bibdiv}
\end{document}